\documentclass{gtmon_a}
\pdfoutput=1

\proceedingstitle{Heegaard splittings of 3-manifolds (Haifa 2005)}
\conferencestart{10 July 2005}
\conferenceend{19 July 2005}
\conferencename{Heegaard splittings of 3-manifolds}
\conferencelocation{Haifa}

\editor{Cameron Gordon}
\givenname{Cameron}
\surname{Gordon}

\editor{Yoav Moriah}
\givenname{Yoav}
\surname{Moriah}

\title{The Heegaard genus of bundles over $S^1$}

\author{Mark Brittenham}
\givenname{Mark}
\surname{Brittenham}
\address{Department of Mathematics\\203 Avery Hall\\\newline
University of Nebraska-Lincoln\\Lincoln, NE 68588-0130\\USA}
\email{mbrittenham2@math.unl.edu}
\urladdr{}

\author{Yo'av Rieck}
\givenname{Yo'av}
\surname{Rieck}
\address{Department of Mathematical Sciences\\University of
Arkansas\\\newline
Fayetteville, AR 72701\\USA}
\email{yoav@uark.edu}
\urladdr{}

\volumenumber{12}
\issuenumber{}
\publicationyear{2007}
\papernumber{02}
\startpage{17}
\endpage{33}

\doi{}
\MR{}
\Zbl{}

\arxivreference{math.GT/0608517}


\received{27 July 2006}
\revised{19 April 2007}
\accepted{20 April 2007}
\published{3 December 2007}
\publishedonline{3 December 2007}
\proposed{}
\seconded{}
\corresponding{}
\version{}

\makeatletter
\def\cnewtheorem#1[#2]#3{\newtheorem{#1}{#3}[section]
\expandafter\let\csname c@#1\endcsname\c@pro}

\let\xysavmatrix\xymatrix
\def\xymatrix{\disablesubscriptcorrection\xysavmatrix}
\AtBeginDocument{\let\bar\wbar\let\tilde\wtilde\let\hat\what}

\newtheorem{pro}{Proposition}[section]
\cnewtheorem{thm}[pro]{Theorem}
\cnewtheorem{lem}[pro]{Lemma}
\cnewtheorem{prop}[pro]{Property}
\cnewtheorem{props}[pro]{Properties}
\cnewtheorem{conds}[pro]{Conditions}
\cnewtheorem{question}[pro]{Question}
\cnewtheorem{example}[pro]{Example}
\cnewtheorem{cnj}[pro]{Conjecture}
\cnewtheorem{rmkk}[pro]{Remark}
\cnewtheorem{rmkks}[pro]{Remarks}
\cnewtheorem{assumptions}[pro]{Assumptions}

\newtheorem{clm}{Claim}

\cnewtheorem{cor}[pro]{Corollary}

\theoremstyle{definition}
\cnewtheorem{dfn}[pro]{Definition}
\cnewtheorem{dfns}[pro]{Definitions}

\theoremstyle{remark}
\newtheorem*{rmk}{Remark}

\makeatother  

\makeautorefname{pro}{Proposition}

\newcommand{\mz}{z_{\mathrm{min}}}
\newcommand{\mm}{\mathrm{min}}

\begin{document}

\begin{htmlabstract}
<p class="noindent">
This paper explores connections between Heegaard genus,
minimal surfaces, and pseudo-Anosov monodromies. Fixing
a pseudo-Anosov map &phi; and an integer n, let M<sub>n</sub>
be the 3&ndash;manifold fibered over S<sup>1</sup> with monodromy &phi;<sup>n</sup>.
</p>
<p class="noindent">
JH Rubinstein showed that for a large enough n every
minimal surface of genus at most h in M<sub>n</sub> is homotopic
into a fiber; as a consequence Rubinstein concludes that
every Heegaard surface of genus at most h for M<sub>n</sub> is
standard, that is, obtained by tubing together two fibers.
We prove this result and also discuss related results of
Lackenby and Souto.
</p>
\end{htmlabstract}

\begin{abstract} 
This paper explores connections between Heegaard genus,
minimal surfaces, and pseudo-Anosov monodromies. Fixing
a pseudo-Anosov map $\phi$ and an integer $n$, let $M_n$
be the 3--manifold fibered over $S^1$ with monodromy $\phi^n$.

JH Rubinstein showed that for a large enough $n$ every
minimal surface of genus at most $h$ in $M_n$ is homotopic
into a fiber; as a consequence Rubinstein concludes that
every Heegaard surface of genus at most $h$ for $M_n$ is
standard, that is, obtained by tubing together two fibers.
We prove this result and also discuss related results of
Lackenby and Souto.
\end{abstract}

\begin{asciiabstract} 
This paper explores connections between Heegaard genus,
minimal surfaces, and pseudo-Anosov monodromies. Fixing
a pseudo-Anosov map phi and an integer n, let M_n
be the 3-manifold fibered over S^1 with monodromy phi^n.

JH Rubinstein showed that for a large enough $n$ every
minimal surface of genus at most h in M_n is homotopic
into a fiber; as a consequence Rubinstein concludes that
every Heegaard surface of genus at most h for M_n is
standard, that is, obtained by tubing together two fibers.
We prove this result and also discuss related results of
Lackenby and Souto.
\end{asciiabstract}

\maketitle

\section{Introduction}
\label{sec:intro}

The purpose of this article to explore theorems of Rubinstein and
Lackenby.  Rubinstein's Theorem studies the Heegaard genus of certain
hyperbolic 3--manifolds that fiber over $S^1$ and Lackenby's Theorem
studies the Heegaard genus of certain Haken manifolds.  Our target
audience is 3--manifold theorists with good understanding of Heegaard
splittings but perhaps little experience with minimal surfaces.  We
will explain the background necessary for these theorems and prove them
(in particular, in \fullref{sec:monotonicity} we explain the main tool
needed for analyzing minimal surfaces).

All  manifolds considered in this paper are closed, orientable 3--manifolds and
all surfaces considered are closed. By the genus of a 3--manifold $M$, denoted
$g(M)$, we mean the genus of a minimal genus Heegaard surface for $M$.  

A {\it least area} surface is a map from a surface into a Riemannian
3--manifold that minimizes the area in its homopoty class.  A {\it minimal
surface} is a critical point of the area functional.  Therefore  a least area
surface is always minimal, as a global minimum is always a critical point.
A local minimum of the area functional is called a stable minimal surface and
has index zero.  However, some minimal surfaces (and in particular the minimal
Heegaard surfaces we will study in this paper) are unstable and have positive
index.  This is similar to a saddle point of the area functional.  An easy
example is the equatorial sphere $\{x_4=0\}$ in $S^3$ (where $S^3$ is the unit
sphere in $\mathbb{R}^4$).  One nice  property that all minimal surfaces share
is that their mean curvature is zero.  This turns out to be equivalent to a
surface being minimal.  It follows that the intrinsic curvature of a minimal
surface is bounded above by the curvature of the ambient manifold.  Thus, the
curvature of a minimal surface $S$ in a hyperbolic manifold is bounded above
by $-1$, and by Gauss--Bonnet the area of $S$ is at most $2\pi\chi(S)$, where
$\chi(S)$ is the Euler characteristic of $S$.  

We assume familiarity with the basic notions of 3--manifold theory
(see, for example, Hempel \cite{hempel} or Jaco \cite{jaco}), the
basic nations about Heegaard splittings (see, for example,
\cite{scharlemann-review}), and Casson and Gordon's concept of {\it
strong irreducibility/weak reducibility} \cite{casson-gordon}.  A more
refined notion, due to Scharlemann and Thompson, is {\it
untelescoping} \cite{untel} (see also Saito, Scharlemann and Schultens
\cite{schschsaito}).  Untelescoping is, in essence, iterated
application of weak reduction (indeed, in some cases a single weak
reduction does not suffice; see Kobayashi \cite{kobayashi-ST}).  In
\fullref{sec:lackenby} we assume familiarity with this concept.

In \cite{rubinstein} Rubinstein used minimal surfaces to study the
Heegaard genus of hyperbolic manifolds that fiber over $S^1$, more
precisely, of closed 3--manifolds that fiber over the circle with
fiber a closed surface of genus $g$ and pseudo-Anosov monodromy (say
$\phi$).  We denote such manifold by $M_\phi$ or simply $M$ when there
is no place for confusion.  While there exist genus two manifolds that
fiber over $S^1$ with fiber of arbitrarily high genus (for example,
consider 0--surgery on 2 bridge knots with fibered exteriors; see
Hatcher and Thurston \cite{hatcher-thurston}) Rubinstein showed that
this is often not the case.  A manifold that fibers over $S^1$ with
genus $g$ fiber has a Heegaard surface of genus $2g+1$ that is
obtained by taking two disjoint fibers and tubing them together once
on each side.  We call this surface and surfaces obtained by
stabilizing it {\it standard}.  $M$ has a cyclic cover of degree $d$
(denoted $M_{\phi^d}$ or simply $M_d$), dual to the fiber, whose
monodromy is $\phi^d$.  Rubinstein shows that for small $h$ and large
$d$ any Heegaard surface for $M_d$ of genus at most $h$ is standard.
In particular, the Heegaard genus of $M_d$ (for sufficiently large
$d$) is $2g+1$.  The precise statement of Rubinstein's Theorem is:

\begin{thm}[Rubinstein]
\label{thm:rubinstein}
Let $M_\phi$ be a closed orientable 3--manifold that fibers over $S^1$ with
pseudo-Anosov monodromy $\phi$.  Let $M_d$ be the $d$--fold cyclic cover of
$M_\phi$ dual to the fiber.

Then for any integer $h \geq 0$ there exists an integer $n>0$ so that for any
$d \geq n$, any Heegaard surface of genus at most $h$ for  $M_d$ is standard.
\end{thm}

\begin{rmk}
In \cite{BS} Bachman and Schleimer gave a combinatorial proof of
\fullref{thm:rubinstein}.  
\end{rmk}

Rubinstein's proof contains two components: the first component is a reduction
to a statement about minimal surfaces.  We state and prove this reduction in
\fullref{sec:reduction}.  It says that if $M_d$ has the property that
every minimal surface of genus at most $h$ is disjoint from some fiber then every
Heegaard surface for $M_d$ of genus at most $h$ is standard.

The second component of Rubinstein's proof is to show that for large enough
$d$, this property holds for $M_d$; this was obtained independently by
Lackenby \cite[Theorem~1.9]{lackenby1}.  A statement and proof are given in
\fullref{sec:main}; we describe it here. 
Let $M$ be a hyperbolic manifold and $F \subset M$ a non-separating surface
(not necessarily a fiber in a fibration over $S^1$).  Construct the $d$--fold
cyclic cover dual to $F$, denoted $M_d$, as follows: let $M^*$ be $M$ cut
open along $F$.  Then $\partial M^*$ has two components, say $F_-$ and $F_+$.
The identification of $F_-$ with $F_+$ in $M$ defines a homeomorphism $h\co F_-
\to F_+$.  We take $d$ copies of $M^*$ (denoted $M^*_i$, with boundaries
denoted $F_{i,-}$ and $F_{i,+}$ ($i=1,\dots,d$))  and glue them together by
identifying $F_{i,+}$ with $F_{i+1,-}$ (the indices
are taken modulo $d$).
The gluing maps are defined using $h$.  The manifold obtained is $M_d$.  In
\fullref{thm:main} we prove that for any $M$   there exists $n$ so that if
$d \geq n$ then any minimal surface of genus at most $h$ in $M_d$ is disjoint
from at least one of the preimages of $F$.  

The proof is an area estimate.
Let $S$ be a minimal surface in a hyperbolic manifold $M_d$ as above;
denote the components of the preimage of $F$ by $F_1,\dots,F_n$.  If $S$
intersects every $F_i$ we
give a lower bound on its area 
by showing that there exists a constant $a > 0$
so that $S$ has area at least $a$ near every $F_i$ that it meets.  Hence if $S$
intersects every $F_i$ it has area at least $ad$.
Fixing $h$, if $d > \frac{2\pi(2h-2)}{a}$ then $S$ has area greater than
$2\pi(2h-2)$.  As mentioned above, the minimal surface $S$ inherits a metric
with curvature bounded above by $-1$, and by Gauss--Bonnet the area of $S$ is
at most $2 \pi (2g(S)-2)$.  Thus $2\pi(2h-2) < \mbox{ area of } (S) \leq 2 \pi
(2g(S)-2)$.  Solving for $g(S)$ we see that $g(S) > h$ as required.  We note
that $a$ is determined by  the geometry of $M$. 

The only tool needed for this is a simple consequence of the {\it Monotonicity
Principle}.  It says that any minimal surface in a hyperbolic ball of radius $R$
that intersects the center of the ball has at least as much area as a
hyperbolic disk of radius $R$.  We briefly explain this in
\fullref{sec:monotonicity}.  For the purpose of illustration we give two
proofs in the case that the minimal surface is a disk.  One of the proofs
requires the following fact: the length of a curve on a sphere or radius $r$
that intersects 
every great circle is at least $2\pi r$, that is, such a curve cannot be shorter
than a great circle.  We give two proofs of this fact in
Appendices~\ref{sec:cruves-on-spheres-1} and \ref{sec:cruves-on-spheres-2}.

Let $N_1$ and $N_2$ be simple manifolds with $\partial N_1 \cong \partial
N_2$ a connected surface of genus $g \geq 2$ (denoted $S_g$).  We emphasize
that by $\partial N_1 \cong \partial N_2$ we only mean that the surfaces are
homeomorphic.   

Let $M'$ be a manifold obtained by gluing $N_1$ to $N_2$ along the boundary.
Then the image of $\partial N_1 = \partial N_2$ (denoted $S$) in $M'$ is an
essential surface.  If $F \subset M'$ is any essential surface with $\chi(F) \geq 0$,
then after isotoping $F$ to minimize $|F \cap S|$, any component of $F
\cap N_1$ or $F \cap N_2$ is essential and has non-negative Euler
characteristic (possibly, $F \cap S = \emptyset$).  But simplicity of $N_1$
and $N_2$ implies that there are no such surfaces.  We conclude that $M'$ is a
Haken manifold with no essential surfaces of non-negative Euler
characteristic.  By Thurston's Uniformization of Haken Manifolds $M'$ is
hyperbolic or Seifert fibered.  If $M'$ is Seifert fibered then $S$ can be
isotoped to be either vertical (that is, everywhere tangent to the fibers) or
horizontal (that is, everywhere transverse to the fibers).  Both cases
contradict simplicity of $N_1$ and $N_2$; the details are left to the
reader.  We conclude that $M'$ is hyperbolic.  Note however, that although
$N_1$ and $N_2$ admit hyperbolic metrics, the
restriction of the hyperbolic metric on $M'$ to $N_1$ and $N_2$ does not have
to resemble them.

After fixing parameterizations $i_1\co S_g \to \partial N_1$ and
$i_2\co S_g \to \partial N_2$ any gluing between  $\partial N_1$ and $\partial
N_2$ is given by a map $i_2 \circ f \circ (i_1^{-1})$ for some map $f\co S_g \to
S_g$.

Fix $f\co S_g \to S_g$ a pseudo-Anosov map, let $M_f$ be the bundle
over $S^1$ with fiber $S_g$ and monodromy $f$, and $M_\infty$ the
infinite cyclic cover of $M_f$ dual to the fiber.  For $n
\in\mathbb{N}$, let $M_n$ be the manifold obtained by gluing $N_1$ to
$N_2$ using the map $i_2 \circ f^n \circ (i_1^{-1})$.  (Note that this
is {\it not} $M_d$.)  Soma \cite{soma} showed that for properly chosen
points $x_n \in M_n$, $(M_n,x_n)$ converge geometrically (in the
Hausdorff--Gromov sense) to $M_\infty$.  In \cite{lackenby} Lackenby
uses an area argument to show that for fixed $h$ and sufficiently
large $n$ every minimal surface of genus at most $h$ in $M_n$ is
disjoint from the image of $\partial N_1 = \partial N_2$ (denoted
$S$).  This implies that any Heegaard surface of genus at most $h$
weakly reduces to $S$, and in particular for sufficiently large $n$,
by Schultens \cite{schultens} $g(M_n) = g(N_1) + g(N_2) - g(S)$.  In
\fullref{sec:lackenby} we discuss Lackenby's Theorem, following the
same philosophy we used for \fullref{thm:rubinstein}.  Finally we
mention Souto's far reaching generalization of Lackenby's Theorem
\cite{souto} and a related theorem of Namazi and Souto
\cite{namazi-souto}; however, a detailed discussion and the proofs of
these theorems are beyond the scope of this note.

\noindent {\bf Acknowledgment}\qua  We thank Hyam Rubinstein for helpful
conversations and the anonymous referee for many helpful suggestions.
 
\section{Reduction to minimal surfaces}
\label{sec:reduction}

In this section we reduce \fullref{thm:rubinstein} to a statement about
minimal surfaces in $M_d$.  We note that the result here applies to any
hyperbolic bundle $M$, but for consistency with applications below we use the
notation $M_d$.

\begin{thm}[Rubinstein]
\label{thm:reduction}
Let $M_d$ be a hyperbolic bundle over $S^1$.  Assume that every minimal
surface of Euler characteristic $\geq 2-2h$ in $M_d$ is disjoint from some
fiber.

Then any Heegaard surface for $M_d$ of genus at most $h$ is standard.
\end{thm}

\begin{proof}
Let $\Sigma \subset M_d$ be a Heegaard surface of genus at most $h$.  By
destabilizing $\Sigma$ if necessary we may assume $\Sigma$ is not stabilized.

Assume first that $\Sigma$ is strongly irreducible.  Then by Pitts and
Rubinstein \cite{pitts-rubinstein} (see also Colding and De Lellis
\cite{colding-lellis}) one of the following holds:
\begin{enumerate}
\item $\Sigma$ is isotopic to a minimal surface.
\item $M_d$ contains a one-sided, non-orientable, incompressible surface (say
  $H$).  Let $H^*$ denote $H$ with an open disk removed.  Then $\Sigma$ is
  isotopic to $\partial N(H^*)$.  Equivalently, $\Sigma$ is isotopic to the
  surface obtained by tubing $\partial N(H)$ once, inside $N(H)$, via a
  straight tube.
\end{enumerate}
Both cases lead to a contradiction:
\begin{enumerate}
\item  Isotope $\Sigma$ to a minimal representative.   Let $\gamma \subset
  M_d$ be a curve.  Since $\Sigma \subset M_d$ is a Heegaard surface $\gamma$
  is freely homotopic into  $\Sigma$.  By assumption, $\Sigma$ is disjoint
  from some fiber $F$.  Thus after free homotopy $\gamma \cap F = \emptyset$,
  and in particular $\gamma$ has algebraic intersection zero with $F$.  But
  this is absurd: clearly there exist a curve $\gamma$ that intersects $F$
  algebraically once.
\item Similarly, any curve $\gamma \subset M_d$ is isotopic into
  $\partial N(H^*)$.  Since $\partial N(H^*) \subset N(H)$ and $N(H)$
  is an $I$--bundle over $H$, $\gamma$ is isotopic into $H$.  Since
  $H$ is essential, by Schoen and Yau \cite{schoen-yau} (see also
  Freedman, Hass and Scott \cite{freedman-hass-scott}) $H$ can be isotoped
  to be least area and in particular minimal.

Note that $2(\chi(H) -1) = 2\chi(H^*) = 2\chi(N(H^*))= \chi(\partial N(H^*)) =
\chi(\Sigma) = 2-2h$.  Hence $\chi(H) = 2-h > 2-2h$.  By assumption $H$ is
disjoint from some fiber $F$.  Thus $\gamma$ can be homotoped to be disjoint
from $F$, contradiction as above.
\end{enumerate}

\begin{rmk}
It is crucial to our proof that $H$ is essential.  Let $H \subset M_d$ be a
non-separating surface so that $\mbox{cl}(M_d  \setminus N(H))$ is a
handlebody.  Let $H^*$ be $H$ with $n$ disks removed, for some $n\geq 1$.  It
is easy to see that $\partial N(H^*)$ is a Heegaard splitting.  However, if
$H$ is compressible, or if $n>1$, then $\partial N(H^*)$ destabilizes.  (The
details are left to the reader.)  The converse was recently studied by
Bartolini and Rubinstein   \cite{bartrubin}.
\end{rmk}

Next assume that $\Sigma$ is weakly reducible.
By Casson and Gordon \cite{casson-gordon} a carefully
chosen weak reduction of $\Sigma$ yields a (perhaps disconnected) essential
surface $S$, and every component of $S$ has genus less than $g(\Sigma)$ 
(and hence less than $h$).  By \cite{schoen-yau} (see also
\cite{freedman-hass-scott}) $S$ is homotopic to a least area (and hence
minimal) representative.  By assumption $S$ is disjoint from some fiber, and
in particular $S$ is embedded in fiber cross $[0,1]$.  Hence $S$ is itself a 
collection of (say $n$) fibers and $\Sigma$ is obtained from $S$ by tubing.  

Note that since $\Sigma$ separates so does $S$.  We conclude that $n$ is
even.  Denote the components of $S$ by $S_1,\dots,S_n$ and the components of
$M_d$ cut open along $S$ by $C_i (i=1,\dots,n$) so that $\partial C_i = F_i
\sqcup F_{i+1}$ (indices taken mod $n$).  Thus $C_i$ is homeomorphic to fiber
cross $[0,1]$.  Fix $i$ and let $\Sigma_i$ be the surface 
obtained by pushing $\partial C_i$ slightly into $C_i$ and then tubing along
the tubes that are contained in $C_i$.  It is easy to see that the component
of $C_i$ cut open along $\Sigma_i$ that contains $\partial C_i$ is a compression
body.  The other component is homeomorphic to a
component obtained by compressing one the handlebodies of $M_d$ cut open
along $\Sigma$.  Hence it is a handlebody.  We conclude that $\Sigma_i$ is a
Heegaard splitting of $C_i$, and both components of $\partial C_i$ are on the
same side 
of $\Sigma_i$.  Scharlemann and Thompson \cite{ST-HS-of-CB} call $\Sigma_i$ a
{\it type II} Heegaard
splitting of $C_i$.  By  \cite{ST-HS-of-CB} either $\Sigma_i$ is obtained by a
single tube that is of the form $\{p\} \times [0,1]$ (for some $p$ in the
fiber) or it is stabilized.  Clearly, if $\Sigma_i$ is stabilized so is $\Sigma$.
We conclude that $\Sigma$ is obtained from $S$ by a single, straight
tube in each $C_i$.

We complete the proof by showing that $n=2$.  Suppose, for a contradiction,
that $n > 2$.
On $S_1$ we see two disks, say $D_0$ and $D_1$, where the tubes in 
$C_0$ and $C_1$ intersect it.  Let $F_1^*$ be $F_1 \setminus (\mbox{int} D_0
\sqcup \mbox{int} D_1)$.  For $i=0,1$ let $\alpha_i \subset F_i^*$ be a
properly embedded 
arc with $\partial \alpha_i \subset \partial D_i$ and so that $|\alpha_0 \cap
\alpha_1| =1$.  Note that $\alpha_i \times [0,1]$ is a meridional disk in $C_i$
($i=0,1$) and these disks intersect once on $F_1$.   Since $n > 2$ these
disks do not have another intersection.  Hence $\Sigma$ destabilizes,
contradicting our assumption.  We
conclude that $n=2$.
\end{proof}

\section{The Monotonicity Principle}
\label{sec:monotonicity}

The Monotonicity Principle studies the growth rate of minimal
surfaces.  All we need is the simple consequence of the Monotonicity
Principle, \fullref{pro:monotonicity}, stated below.  For illustration
purposes, we give two proofs of \fullref{pro:monotonicity} in the
special case when the minimal surface intersects the ball in a
(topological) disk.  A proof for the Monotonicity Principle for annuli
is given in Lackenby \cite[Section~6]{lackenby1}.  For the general
case, see Simon \cite{simon} or Choe \cite{choe}.

We will use the following facts about minimal surfaces: (1) if a minimal surface
$F$ intersects a small totally geodesic disk $D$ and locally  $F$ is contained
on one side of $D$ then $D \subseteq F$.  (2) If $D$ is a little piece of the round
sphere $\partial B$ (for some metric ball $B$) and $F \cap D \neq \emptyset$,
then locally $F \not\subset B$.  Roughly speaking, these facts state that a
minimal surface cannot have ``maxima'' (or, the maximum principle for minimal
surfaces).

In this section we use the following notation: $B(r)$ is a hyperbolic ball of
radius $r$, which for convenience we identify with the ball of radius $r$ in
the Poincar\'e ball model in $\mathbb{R}^3$, centered at $O = (0,0,0)$.   The
boundary of 
$B(r)$ is denoted $\partial B(r)$.  A great circle in $\partial B(r)$ is the
intersection of $B(r)$ with a totally geodesic disk that
contains $O$, or, equivalently, the intersection of $\partial B(r)$ with a
2--dimensional subspace of $\mathbb{R}^3$.  For
convenience, we use the horizontal circle (which we shall call the equator) as a
great circle and denote the totally geodesic disk it bounds $D_0$.
Note that $\partial D_0$ separates $\partial B(r)$ into two disks which we
shall call the northern and
southern hemispheres, and $D_0$ separates $B(r)$ into two (topological) balls
which we shall call the northern and
southern half balls.  The ball $B(r)$ is foliated by geodesic disks
$D_t$ ($-r \leq t \leq +r$), where $D_t$ is the intersection of $B(r)$ with the
geodesic plane that is perpendicular to the $z$--axis and intersects it at
$(0,0,t)$.  Here and throughout this paper, we denote the area of a hyperbolic
disk of radius $r$ by $a(r)$.  In the first proof below we use the fact that
if a curve on a sphere intersects every great circle then it is at least as
long as a great circle (\fullref{pro:great-circles}).  This is an
elementary fact in spherical geometry.  In
Appendices~\ref{sec:cruves-on-spheres-1} and \ref{sec:cruves-on-spheres-2} we
give two proofs of this fact, however, we encourage the reader to find her/his
own proof and send it to us.

\begin{pro} 
\label{pro:monotonicity} 
Let $B(R)$ be a hyperbolic ball of radius $R $ centered at $O$ and $F \subset
M$ a minimal surface so that $O \in F$.  Then the area of $F$ is at least
$a(R)$. 
\end{pro} 

\begin{rmk}
Lackenby's approach \cite{lackenby1} does not require the full strength of the
Monotonicity Principle.  He only needs the statement for annuli, and in that
case he gives a self-contained proof in  \cite[Section~6]{lackenby1}.
\end{rmk}

We refer the reader to \cite{simon} or \cite{choe} for a proof.  For the
remainder of the section, assume $F \cap B(R)$ is topologically a disk.  Then
we have: 

\begin{proof}[First proof]
Fix $r$, $0 < r \leq R$. Fix a great circle in $\partial B(r)$ (which
for convenience we identify with the equator).  Suppose that $F \cap \partial
B(r)$ is not the equator, we will show that $F \cap \partial B(r)$ intersects
both the northern and southern hemispheres.  Suppose for contradiction for
some $r$ this is not the case.  Then one of the following holds: 

\begin{enumerate}
\item $F \cap \partial B(r) = \emptyset$.
\item $F \cap \partial B(r) \neq \emptyset$ and $F \cap \partial B(r)$ does
  not intersect one of the two hemispheres. 
\end{enumerate}

Assuming Case~(1) happens, and let $r' > 0$ be the largest value for which $F
\cap \partial B(r') \neq \emptyset$.  Then $F$ and $\partial B(r')$ contradict
fact (2) mentioned above.

Next assume Case~(2) happens (say $F$ does not intersect the southern
hemisphere).  Let $t$  be the  most negative value for which $F \cap D_t \neq
\emptyset$.  Since $O \in F$, $-r < t \leq 0$.  Then by fact~(1) above, $F$
must coincide with $D_t$.  If $t < 0$ then $D_t$ intersects the southern
hemisphere, contrary to our assumptions.  Hence $t=0$ and $F$ is itself $D_0$;
thus $F \cap B(r)$ is the equator, again contradicting our assumptions.

By assumption $F \cap B(R)$ is a disk and therefore $F \cap \partial B(r)$ is
a circle.  Clearly, a circle that intersects both the northern and the
southern hemispheres must intersect the equator.  We conclude that $F \cap
\partial B(r)$ intersects the equator, and as the equator was chosen
arbitrarily, $F \cap \partial B(r)$ intersects every great circle.  By
\fullref{pro:great-circles} $F \cap B(r)$ is at least as long as 
a great circle in  $\partial B(r)$.   Since the intersection of a totally
geodesic disk with $\partial B(r)$ is a great circle,
integrating these lengths shows that
the area of $F \cap B(r)$ grows at least as fast as the area of a  geodesic
disk, proving the proposition.
\end{proof}

\begin{proof}[Second proof] 
Restricting the metric from $M$ to $F$,  distances can increase but cannot
decrease.  Therefore $F \cap \partial B(R)$ is at distance (on $F$) at least
$R$ from  $O$ and we conclude that $F$ contains an entire disk of radius $R$.
The induced metric on $F$ has curvature at most $-1$ and therefore areas on $F$
are at least as big as areas in $\mathbb{H}^2$.  In particular, the disk of
radius $R$ about $O$ has area at least $a(R)$.
\end{proof}

\section{Main Theorem}
\label{sec:main}

By Theorem~\ref {thm:reduction} the main task in proving
\fullref{thm:rubinstein} is showing that (for large enough $d$) a minimal
surface of genus at most $h$ in $M_d$ is disjoint from some fiber $F$.  Here
we prove: 

\begin{thm}
\label{thm:main}
Let $M$ be a compact, orientable hyperbolic manifold and $F \subset M$ a
non-separating, orientable surface.
Let $M_d$ denote the cyclic cover of $M$ dual to $F$ of degree $d$ (as in the
introduction).

Then for any integer $h \geq 0$ there exists a constant $n$ so that for $d
\geq n$, any minimal surface of genus at most $h$ in $M_d$ is disjoint from a
component of the preimage of $F$. 
\end{thm}

\begin{proof}
Fix an integer $h$.

Denote the distance in $M$ by $d(\cdot,\cdot)$.  Push $F$ off itself to obtain
$\widehat{F}$, a surface parallel to $F$ and disjoint from it.  For each point
$p \in F$ define:
$$R(p) = \min \{\mbox{radius of injectivity at }p, \ d(p,\widehat{F}) \}.$$%
Since $\widehat{F}$ is compact $R(p) > 0$.  Define:
$$R = \min \{R(p) |p \in F\}.$$  %
Since $F$ is compact $R > 0$.  Note that $R$ has the following property: for
any $p \in F$, the set $\{q \in M: d(p,q) < R \}$ is an embedded ball and this
ball is disjoint from $\widehat{F}$.  As above, let $a(R)$ denote the area of
a hyperbolic disk of radius $R$.

Let $n$ be the smallest integer bigger than $\frac{2\pi (2h-2)}{a(R)}$.  Fix
an integer $d \geq n$.  Denote the preimages of $F$ in $M_d$ by
$F_1,\dots,F_d$.

Let $S$ be a minimal surface in $M_d$.  Suppose $S$ cannot be isotoped to be
disjoint from the preimages of $F_i$ for any $i$.  We will show that $g(S) >
h$, proving the theorem.

Pick a point $p_i \in F_i \cap S$ ($i=1,\dots,d$) and let $B_i$ be the set
$\{p \in M_d | d(p,p_i) < R\}$.  By choice of $R$, for each $i$, $B_i$ is an
embedded ball and the preimages of $\widehat{F}$ separate these balls; hence
for $i \neq j$ we see that $B_i \cap B_j = \emptyset$.  $S
\cap B_i$ is a minimal surface in $B_i$ that intersects its center and by
\fullref{pro:monotonicity} (the Monotonicity Principle) has area at
least $a(R)$.  Summing these areas we see that the area of $S$ fulfills:
\begin{eqnarray*}
         \mbox{Area of } S &\geq&   d \cdot a(R) \\ &\geq& n \cdot a(R) \\ &>&
              \frac{2\pi (2h-2)}{a(R)} \cdot a(R) \\ &=& 2 \pi (2h-2)
\end{eqnarray*}
But a minimal surface in a hyperbolic manifold has curvature $\leq -1$ and
hence by the Gauss--Bonnet Theorem, the area of $S \leq -2 \pi \chi(S) =
2\pi(2g(S - 2))$. Hence, the genus of $S$ is greater than $h$.
\end{proof}

\begin{rmkk}[Suggested project]{\rm
In \fullref{thm:main} we treat the covers dual to a non-separating
essential surface (denoted $M_d$ there).  
In the section titled ``Generalization'' of \cite{lackenby}, Lackenby
shows (among other things) how to amalgamate along non-separating surfaces.
Does his construction and \fullref{thm:main} give useful bounds on the
genus of $M_d$, analogous to \fullref{thm:rubinstein}?   }\end{rmkk}

\section{Lackenby's Theorem}
\label{sec:lackenby}

Lackenby studied the Heegaard genus of manifolds containing separating
essential surfaces.  Here too, the result is asymptotic.  We begin by
explaining the set up.  Let $N_1$ and $N_2$ be simple manifolds with
$\partial N_1 \cong \partial N_2$ a connected surface of genus $g \geq 2$
(that is, $\partial N_1$ and $\partial N_2$ are homeomorphic).
Let $S$ be a surface of genus $g$ and $\psi_i\co S \to \partial N_i$ 
parameterizations of the boundaries ($i=1,2$).  
Let $f\co S \to S$ be a pseudo-Anosov map.
For any $n$ we construct the map $f_n = \psi_2 \circ f^n \circ (\psi_1)^{-1}:
\partial N_1 \to \partial N_2$.  By identifying $\partial N_1$ with $\partial
N_2$ by the map $f_n$ we obtain a closed hyperbolic manifold $M_n$.   Let $S
\subset M_n$ be the image of $\partial N_1 =\partial N_2$.  With this we are
ready to state Lackenby's Theorem:

\begin{thm}[Lackenby \cite{lackenby}]
\label{thm:lackenby}
With notation as in the previous paragraph, for any $h$ there exists $N$ so
that for any $n \geq N$ any genus $h$ Heegaard surface for $M_n$ weakly
reduces to $S$.  In particular, by setting $h = g(N_1) + g(N_2) - g(S)$ we see
that there exists $N$ so that if $n \geq N$ then $g(M_n) = g(N_1) + g(N_2) -
g(S)$.
\end{thm}

\begin{proof}[Sketch of proof]
As in Sections~\ref{sec:reduction} and \ref{sec:main}, the proof has two parts
which we bring here as two claims:

\begin{clm}
Suppose that every every minimal surface in $M_n$ of genus at most $h$ can be
homotoped to be disjoint from $S$.  Then any Heegaard surface of genus at most
$h$ weakly reduces to $S$.  In particular, if $h \geq g(N_1) + g(N_2) - g(S)$
then $g(M_n) = g(N_1) + g(N_2) - g(S)$.
\end{clm} 

\begin{clm}
There exists $N$ so that if $n \geq N$ then any minimal surface of genus at
most $h$ in $M_n$ can be homotoped to be disjoint from $S$.
\end{clm}

Clearly, Claim~1 and 2 imply Lackenby's Theorem.  We now sketch their proofs.

We paraphrase Lackenby's proof of Claim~1: let $\Sigma$ be a Heegaard surface
of genus at most $h$.  Then by Scharlemann and Thompson \cite{st-untel} $\Sigma$
untelescopes to a collection of connected surfaces $F_i$ and $\Sigma_j$ where
$\cup_i F_i$ is an essential surface (with $F_i$ its components) and
$\Sigma_j$ are strongly irreducible Heegaard surfaces for the components of
$M_n$ cut open along $\cup_i F_i$; in particular $M_n$ cut open along $(\cup_i
F_i) \cup (\cup_i \Sigma_j)$ consists of compression bodies and the images of
the $F_i$'s form $\partial_-$ of these compression bodies.  Since $F_i$
and $\Sigma_j$ are obtained by compressing $\Sigma$, they all have genus less
than $h$.  

By \cite{schoen-yau}, \cite{freedman-hass-scott}, and \cite{pitts-rubinstein}
the surfaces $F_i$ and $\Sigma_j$ can be made minimal.  We explain this
process here: since $F_i$ are essential surfaces they can be made minimal by
\cite{schoen-yau} (see also \cite{freedman-hass-scott}).   Next, since the
$\Sigma_j$'s are strongly irreducible Heegaard surfaces for the components of
$M_n$ cut open along $\cup_i F_i$, each $\Sigma_j$ can be made minimal within its
component by \cite{pitts-rubinstein} (see also \cite{colding-lellis}).   Note
that the surfaces $F_i$ and $\Sigma_j$ are disjointly embedded.

By assumption, $S$ can be isotoped to be disjoint from every $F_i$ and
every $\Sigma_j$.  Therefore, $S$ is an essential closed surface in a
compression body and must be parallel to a component of
$\partial_-$. Therefore, for some $i$, $S$ is isotopic to $F_i$.  In
Rieck and Kobayashi \cite[Proposition~2.13]{rieck-kobayashi} it was
shown that if $\Sigma$ untelescopes to the essential surface $\cup_i
F_i$, then $\Sigma$ weakly reduces to any {\it connected separating}
component of $\cup_i F_i$; therefore $\Sigma$ weakly reduces to $S$.
This proves the first part of Claim~1.

Since $S$ is connected any minimal genus Heegaard splittings for $N_1$ and
$N_2$ can be amalgamated (the converse of weak reduction \cite{schultens}).
By amalgamating minimal genus Heegaard surfaces we see  that for any $n$,
$g(M_n) \leq g(N_1) + g(N_2) - g(S)$.  By applying the first part of Claim~1
with $h = g(N_1) + g(N_2) - g(S)$ we see that for sufficiently large $n$,
a minimal genus Heegaard surface for $M_n$ weakly reduces to $S$; by
\cite[Proposition~2.8]{rieck-kobayashi}
$g(M_n)=  g(N_1) + g(N_2) - g(S)$, completing the proof of Claim~1.

We now sketch the proof of Claim~2.  Fix $h$ and assume that for arbitrarily
high values of $n$, $M_n$ contains a minimal surface (say $P_n$) of genus
$g(P_n) \leq h$ that  cannot be homotoped to be disjoint from $S$.  Let $M_f$
be the bundle over $S^1$ with monodromy $f$ and fix two disjoint fibers $F$,
$\widehat{F} \subset M_f$.  Let $R$ be as in \fullref{sec:main}.  Let
$M_{\infty}$ be the infinite cyclic cover dual to the fiber. Soma \cite{soma}
showed that there are points $x_n \in M_n$ so that $(M_n,x_n)$ converges in
the sense of Hausdorff--Gromov to the manifold $M_\infty$.  These points are
near the minimal surface $S$, and the picture is that $M_n$ has a very long
``neck'' that looks more and more like $M_\infty$.

For sufficiently large $n$ there is a ball $B(r) \subset M_n$ for arbitrarily
large $r$ that is $1-\epsilon$ isometric to $B_\infty(r) \subset M_\infty$.
Note that $B_\infty(r)$ contains arbitrarily many lifts of $F$ separated by
lifts of $\widehat{F}$.  Since $P_n$ cannot be isotoped to be disjoint from
$S$, its image  in $M_\infty$ cannot be isotoped off the preimages of $F$.  As
in \fullref{sec:main} we conclude that the images of $P_n$ have
arbitrarily high area.   However, areas cannot be distorted arbitrarily by a
map that is $1-\epsilon$ close to an isometry.  Hence the areas of $P_n$ are
unbounded, contradicting Gauss--Bonnet; this contradiction completes our
sketch.
\end{proof}

In \cite{souto} Souto generalized Lackenby's result (see also a recent
paper by Li \cite{li}).  Although his work is beyond the scope of this
paper, we give a brief description of it here.  Instead of powers of
maps, Souto used a combinatorial condition on the gluings: fixing
essential curves $\alpha_i \subset N_i$ ($i=1,2$) and $h>0$, Souto
shows that if $\phi\co N_1 \to N_2$ fulfills the condition
``$d_{\mathcal{C}}(\phi(\alpha_1), \alpha_2)$ is sufficiently large''
then any Heegaard splitting for $N_1 \cup_\phi N_2$ of genus at most
$h$ weakly reduces to $S$.  The distance Souto
uses---$d_{\mathcal{C}}$---is the distance in the ``curve complex''
(as defined by Hempel \cite{hempel-distance}) and {\it not} the
hyperbolic distance.  Following Kobayashi \cite{kobayashi-height}
Hempel showed that raising a fixed monodromy $\phi$ to a sufficiently
high power does imply Souto's condition.  Hence Souto's condition is
indeed weaker than Lackenby's, and it is in fact too weak for us to
expect Soma-type convergence to $M_\infty$.  However, using Minsky
\cite{minsky} Souto shows that given a sequence of manifolds
$M_{\phi_n}$ with $d_{\mathcal{C}}(\phi_n(\alpha_1), \alpha_2) \to
\infty$, the manifolds $M_{\phi_n}$ are ``torn apart'' and the cores
of $N_1$ and $N_2$ become arbitrarily far apart.  For a precise
statement \cite[Proposition~6]{souto}.  Souto concludes that for
sufficiently large $n$, any minimal surface for $M_n$ that intersects
both $N_1$ and $N_2$ has high area and therefore genus greater than
$h$.  Souto's Theorem now follows from Claim~1 above.

A similar result was obtained by Namazi and Souto \cite{namazi-souto}
for gluing of handlebodies.  They show that if $N_1$ and $N_2$ are genus $g$
handlebodies and $\partial N_1 \to \partial N_2$ is a generic pseudo-Anosov
map (for a precise definition of ``generic'' in this case see
\cite{namazi-souto}) then for any  $\epsilon > 0$ and for large enough $n$ the
manifold $M_{f^n}$ obtained by gluing $N_1$ to $N_2$ via $f^n$ admits a
negatively curved metric with curvatures $K$ so that $-1 - \epsilon < K < -1 +
\epsilon$.  Namazi and Souto use this metric to conclude many things about
$M_{f^n}$, for example, that both its Heegaard genus and its rank (that is,
number of generators needed for $\pi_1(M_{f^n})$) are exactly $g$.

\appendix

\section{Appendix: Short curves on round spheres: take one}
\label{sec:cruves-on-spheres-1}

In this section we prove the following proposition, which is a simple exercise
in spherical geometry used in \fullref{sec:monotonicity}.  Let $S^2(r)$ be a
sphere of constant curvature $+(\frac{1}{r})^2$.  We isometrically identify
$S^2(r)$ with $\{(x,y,z) \in \mathbb{R}^3 | x^2 + y^2 + z^2 = r^2\}$ and refer
to it as a round sphere of radius $r$.

\begin{pro}
\label{pro:great-circles}
Let $S^2(r)$ be a round sphere of radius $r$ and $\gamma \subset S^2$ a rectifiable
closed curve. Suppose $l(\gamma) \leq 2\pi r$ (the length of great circles). Then
$\gamma$ is disjoint from some great circle.  
\end{pro} 

\begin{rmk}
The proof also shows that if $\gamma$ is a {\it smooth} curve that meets every
great circle then $l(\gamma) = 2\pi r$ if and only if $\gamma$ is itself a
great circle.
\end{rmk}

\begin{proof}
Let $\gamma$ be a curve that intersects every great circle.  Let $\mz$ (for
some $\mz \in \mathbb{R}$) be the minimal value of the $z$--coordinate, taken over
$\gamma$.  Rotate $S^2(r)$ to maximize $\mz$.  If $\mz > 0$ then $\gamma$ is
disjoint from the equator, contradicting our assumption.  We assume from now
on $\mz \leq 0$.

Suppose first $\mz = 0$.  Suppose, for contradiction, that there exists a
closed arc $\alpha$ on the equator so that $l(\alpha) = \pi r$ and $\alpha
\cap \gamma = \emptyset$.  By rotating $S^2(r)$ about the $z$--axis (if
necessary) we may assume $\alpha  = \{(x,y,0) \in S^2(r) | y \leq 0\}$.  Then
rotating $S^2(r)$ slightly about the $x$--axis pushes the points  $\{(x,y,0)
\in S^2(r) | y > 0\}$ above the $xy$--plane.  By compactness of $\gamma$ and
$\alpha$ there is some $\epsilon$ so that $d(\gamma,\alpha) > \epsilon$.
Hence if the rotation is small enough, no point of $\gamma$ is moved to (or
below) $\alpha$.  Thus, after rotating $S^2(r)$, $\mz > 0$, contradiction.  We
conclude that every arc of the equator of length $\pi r$ contains a point of
$\gamma$.  Therefore there exists a sequence of points $p_i \in \gamma \cap
\{(x,y,0)\}$ ($i=1,\dots,n$, for some $n \geq 2$), ordered by their  order
along the equator ({\it not}  along $\gamma$), so that $d(p_i,p_{i+1})$ is at
most half the equator (indices taken modulo $n$).  The shortest path
connecting $p_i$ to $p_{i+1}$ is an arc of the equator, and we conclude that
$l(\gamma) \geq 2\pi r$ as required.  If we assume, in addition, that
$l(\gamma) = 2\pi r$ then either $\gamma$ is itself the equator or   $\gamma$
consists of two arcs of great circle meeting at $c_1 \cup c_2$.  Note that
this can in fact happen, but then $\gamma$ is not smooth.  This completes the
proof in the case $\mz = 0$

Assume next $\mz < 0$.  Let $c_{\mm}$ be the latitude of $S^2(r)$ at $z =
\mz$,  and denote the length of $c_{\mm}$ by $d_{\mm}$.  Suppose there is an
open arc of $c_{\mm}$ of length $\frac{1}{2} d_{\mm}$ that does not intersect
$\gamma$.  Similar to above, by rotating $S^2(r)$ we may assume this arc is
given by $\{(x,y,\mz) \in c_{\mm} | y < 0\}$.  Then a tiny rotation about the
$x$--axis increases the $z$--coordinate of all points $\{(x,y,z)| y \geq 0, \ z
\leq 0\}$.  As above ,this increases $\mz$, contradicting our choice of $\mz$.
Therefore there is a collection of points $p_i \in \gamma \cap c_{\mm}$
($i=1,\dots,n$, for some $n \geq 3$), ordered by their order along the equator
({\it not}  along $c_{\mm}$), so that $d(p_i,p_{i+1}) < \frac{1}{2} d_{\mm}$
(indices taken modulo $n$).  The shortest path connecting $p_i$ to $p_{i+1}$
is an arc of a great circle.  However, such arc has points with $z$--coordinate
less than $\mz$, and therefore cannot be a part of $\gamma$.  The shortest
path containing all the $p_i$'s on the punctured sphere  on $\{(x,y,z) \in
S^2(r) | z \geq \mz\}$ is the boundary, that is, $c_{\mm}$ itself.
Unfortunately, $l(c_{\mm}) < 2\pi r$.  Upper hemisphere to the rescue!
$\gamma$ must have a point with $z$--coordinate at least $-\mz$, for otherwise
rotating $S^2(r)$ by $\pi$ about any horizontal axis would decrease $\mz$.
Then $l(\gamma)$ is at least as long as the shortest curve containing the
$p_i$'s and some point $p$ on or above $c_{\mm}$, the circle of $\gamma$ at $z
= \mz$.  Let $\gamma$ be such a curve.  By reordering the indices if necessary
it is convenient to assume that $p$ is between $p_1$ and $p_2$.  It is clear
that moving $p$ so that its longitude is between the longitudes of $p_1$ and
$p_2$ shortens $\gamma$ (note that since $d(p_1,p_2) < \frac{1}{2} d_{\mm}$
this is well-defined).  We now see that $\gamma$ intersects the equator in two
point, say $x_1$ and $x_2$.  Replacing the two arcs of $\gamma$ above the
equator by the short arc of the equator decreases length.  It is not hard to
see that the same hold when we replace the arc of $\gamma$ below the equator
with the long arc of the equator.  We conclude that $l(\gamma) >
l(\mbox{equator}) = 2\pi r$.
\end{proof}

\section{Appendix: Short curves on round spheres: take two}
\label{sec:cruves-on-spheres-2}

We now give a second proof of \fullref{pro:great-circles}.  For
convenience of presentation we take $S^2$ to be a sphere of radius 1.  Let
$\gamma$ be a closed curve that intersects every great circle.  Every great
circle is defined by two antipodal points, for example, the equator is defined
by the poles.  Thus, the space of great circles is $\mathbb{R}P^2$.  Since
$S^2$ has area $4\pi$, $\mathbb{R}P^2$ has area $2\pi$.  Let $f\co S^2 \to
\mathbb{R}P^2$ be the ``map'' that assigns to a point $p$ all the great
circles that contain $p$; thus, for example, if $p$ is the north pole then
$f(p)$ is the projection of the equator to $\mathbb{R}P^2$.

Let $C$ be a great circle.  We claim that $\gamma \cap C$ contains at least two
points of $\gamma$.  (If $\gamma$ is not embedded then the two may be the same
point of $C$.)  Suppose, for a contradiction, that $\gamma$ meets some great
circle (say the equator) in one point only (Say $(1,0,0)$).  By the
Jordan Curve Theorem, $\gamma$ does not cross the equator.  
By tilting the equator slightly about the $y$--axis it
is easy to obtain a great circle disjoint from $\gamma$.  Hence we see that
$\gamma$ intersects every great circle at least twice.  Equivalently,
$f(\gamma)$ covers $\mathbb{R}P^2$ at least twice.

Let $\alpha_i$ be a small arc of a great circle, of length
$l(\alpha_i)$; note that this length is exactly the angle $\alpha_i$ supports
in radians.  Say for convenience $\alpha_i$ starts at the north pole and goes
towards the equator.  The points that define great circles that intersect
$\alpha_i$ are given by tilting the equator by $\alpha_i$ radians.  This gives
a set whose area is $\alpha_i/\pi$ of the total area of $S^2$.  Since the area
of $S^2$ is $4\pi$, it gives a set of area $4l(\alpha_i)$.  This set is
invariant under the antipodal map, and so projecting to $\mathbb{R}P^2$ the
area is cut by half, and we get:
\begin{equation}
\label{eq:areas_and_lengths}
\mbox{Area of } f(\alpha_i) = 2l(\alpha_i).
\end{equation}
Fix $\epsilon > 0$.  Let $\alpha$ be an approximation of $\gamma$ by small
arcs of great circles, say $\{\alpha_i\}_{i=1}^n$ are the segments of
$\alpha$.  We require $\alpha$ to approximate $\gamma$ well in the following
two senses:

\begin{enumerate}
\item $l(\alpha) \leq l(\gamma) + \epsilon$.
\item Under $f$, $\alpha$ covers $\mathbb{R}P^2$ as well as $\gamma$ does
  (except, perhaps, for a set of measure $\epsilon$); ie, the area of
  $f(\alpha) \geq $ the area of $f(\gamma) - \epsilon$ (area measured with
  multiplicity).
\end{enumerate}
From this we get:
\begin{eqnarray*}
  4\pi - \epsilon &=& \mbox{twice the area of } \mathbb{R}P^2 - \epsilon \\%
  &\leq& \mbox{the area of } f(\gamma) - \epsilon \\%
       &\leq& \mbox{area of } f(\alpha) \\
       &=& \Sigma_{i=1}^n \mbox{area of }f(\alpha_i) \\%
       &=& \Sigma_{i=1}^n 2l(\alpha_i) \\%
       &=& 2l(\alpha) \\%
       &\leq& 2(l(\gamma) + \epsilon).%
\end{eqnarray*}

(In the fifth equality we use Equation~\eqref{eq:areas_and_lengths}.)  Since
$\epsilon$ was arbitrary, dividing by 2 we get the desired result: $2\pi
\leq l(\gamma)$.

\bibliographystyle{gtart}
\bibliography{link}

\end{document}